\documentclass[12pt]{amsart}
\usepackage{amssymb,latexsym,eufrak,amsmath,amscd,graphicx}
 \usepackage[all]{xy}

\theoremstyle{plain}
\newtheorem{theorem}{Theorem}[section]
\newtheorem{proposition}[theorem]{Proposition}
\newtheorem{corollary}[theorem]{Corollary}
\newtheorem{lemma}[theorem]{Lemma}

\newtheorem*{theoremnn}{Theorem}

\newtheorem{theoremalpha}{Theorem}
\newtheorem{corollaryalpha}[theoremalpha]{Corollary}

\theoremstyle{definition}

\newtheorem{remark}[theorem]{Remark}
\newtheorem{example}[theorem]{Example}

\newcommand{\lra}{\longrightarrow}
\newcommand{\noi}{\noindent}
\newcommand{\AAA}{\mathbb{A}}

\newcommand{\RR}{\mathbb{R}}
\newcommand{\NN}{\mathbb{N}}

\newcommand{\CC}{\mathbb{C}}
\newcommand{\QQ}{\mathbb{Q}}

\newcommand{\OO}{\mathcal{O}}

\newcommand{\II}{\mathcal{I}}

\newcommand{\fra}{\frak{a}}
\newcommand{\frb}{\frak{b}}

\newcommand{\lct}{\textnormal{lct}}
\newcommand{\Spec}{\textnormal{Spec}}

\newcommand{\ord}{\textnormal{ord}}
\newcommand{\codim}{\textnormal{codim}}

\newcommand{\pr}{\prime}

\newcommand{\Cont}{\textnormal{Cont}}

\newcommand{\supp}{\textnormal{supp}}
\newcommand{\val}{\textnormal{val}}
\begin{document}

\title[Contact Loci in Arc Spaces]  {Contact Loci in
Arc Spaces} 

\author[Ein]{Lawrence Ein}
\address{Department of Mathematics \\University of
Illinois at Chicago \hfil\break\indent 
851 South Morgan Street (M/C 249)\\ Chicago, IL 60607-7045, USA}
\email{ein@math.uic.edu}

\author[Lazarsfeld]{Robert Lazarsfeld} 
\address{Department of Mathematics
\\ University of Michigan \\ Ann Arbor, MI 48109, USA}
\email{rlaz@umich.edu}

\author[Musta\c{t}\v{a}]{Mircea~Musta\c{t}\v{a}}\address{Department
of Mathematics
\\ Harvard University \\ One Oxford Street \\ Cambridge, MA 02138, USA}
\email{mirceamustata@yahoo.com}
\maketitle



\section*{Introduction}

The purpose of this paper is to study  the loci of arcs
on a smooth variety defined  by order of contact with a
fixed subscheme. Specifically, we establish a Nash-type
correspondence showing that the irreducible components
of these loci arise from (intersections of)  exceptional
divisors in a resolution of singularities. We show also
that these loci account for all the valuations
determined by irreducible cylinders in the arc space.
 Along the way, 
we recover in an elementary fashion   -- without using 
motivic integration --  results of the third
author from
\cite{Mu1} and
\cite{Mu2} relating singularities to arc spaces. Moreover, we extend these
results to give  a jet-theoretic interpretation of multiplier ideals.

Let $X$ be a smooth complex variety of dimension $d$.
Given $m \ge 0$ we denote by
\[  X_m \ = \ \text{Hom} \big( \,  \Spec \, \CC[t]/
(t^{m+1})
\, , \, X \, \big) \]
the space of $m^{\text{th}}$ order arcs on $X$. Thus
$X_m$ is a smooth variety of dimension $d(m+1)$, and
the truncation morphism  $\tau_{m+1,m} : X_{m+1} \lra
X_m$ realizes each of these  spaces as a $\CC^d$-bundle
over the previous one.  The inverse limit
$X_{\infty} $ of the
$X_m$ parametrizes all formal arcs on
$X$, and one writes $\psi_m : X_{\infty} \lra X_m$ for
the natural map. These spaces have attracted a great
deal of interest in recent years thanks to their central
role in the theory of motivic integration (see, for example, \cite{batyrev},
\cite{denef}, \cite{denef2}). In his papers \cite{Mu1}  and \cite{Mu2}
the third author used this machine to give  
arc-theoretic interpretations  of some of the basic 
invariants of higher dimensional geometry.

Consider now a non-zero ideal sheaf $\fra \subseteq
\OO_X$ defining a subscheme $Y \subseteq X$. Given a
finite or infinite arc $\gamma$ on $X$, the
 \text{order of vanishing} 
 of $\fra$  --- or
the \text{order of contact} of the corresponding
scheme
$Y$ --- along
$\gamma$ is defined in the natural
way.\footnote{Specifically,  pulling $\fra $ back  via
$\gamma$ yields an ideal $(t^e)$ in $\CC[t]/(t^{m+1})$
or
$\CC[[t]]$, and one sets
\[  \ord_\gamma(\fra) \ = \ \ord_\gamma(Y) \ = \ e.
\] (Take  $\ord_\gamma(\fra) = m+1$ when 
$\fra$ pulls back to the zero ideal in
$\CC[t]/(t^{m+1})$ and $\ord_\gamma(\fra) = \infty$
when it pulls back to the zero ideal in $\CC[[t]]$.)}
For a fixed integer  
$p
\ge 0$, we define the \textit{contact loci}
\[ \Cont^p(Y) \ = \ \Cont^p(\fra) \ = \ \big\{ \, \gamma
\in X_{\infty} \, | \, \ord_{\gamma}(\fra) = p \, \big
\}. \]  These are locally closed cylinders, i.e. 
they arise as the common pull-back  of
the locally closed sets
\begin{equation}
 \Cont^p(Y)_m \ = \  \Cont^p(\fra)_m \ =_{\text{def}} \
\big
\{
\,\gamma
\in X_m \mid \ \ord_\gamma(\fra) = p \, \big \}
\end{equation}
defined for any $m \ge p$. 
By an irreducible component $W$ of $\Cont^p(\fra)$ we
mean the inverse image of a component $W_m$ of
$\Cont^p(\fra)_m$ (these being in bijection for all $m
\ge p$ via the truncation maps).  The codimension of
$W$ in
$X_{\infty}$ is  the common codimension of any of the
$W_m$ in $X_m$.

Our first result establishes a Nash-type
correspondence showing that the irreducible components
of
$\Cont^p(\fra)$ arise from exceptional divisors 
appearing in  an embedded resolution of $Y$.  
Specifically, fix a log resolution
$\mu : X^\pr \lra X$ of
$\fra$. Thus
$X^\pr$ is a smooth variety which carries a simple
normal crossing divisor $E = \sum_{i=1}^t E_i  $,
$\mu$ is a projective birational map, and 
\begin{equation}  \fra^\pr \ = \ \fra \cdot \OO_{X^\pr}
\ = \
\OO_{X^\pr} (- \sum_{i=1}^t r_i E_i). \end{equation}
For simplicity, we assume here that all $r_i$ are positive
(we can do this if $\mu$ is an isomorphism over $X\setminus Y$),
and leave the general
case for the main body of the paper.
We write
\begin{equation} \label{Defining.K.Eqn}
K_{X^\pr/X} \ = \ \sum_{i=1}^t \, k_i \; E_i
\end{equation} where $K_{X^\pr/X} = K_{X^\pr} - \mu^*
K_X$ is the relative canonical divisor defined by the
vanishing of
$\det d\mu$.  After perhaps some further blowings-up, we
may -- and do -- assume that any non-empty intersection
of the
$E_i$ is connected and hence irreducible. Observe that
$\mu$ gives rise in the natural way to morphisms
$
\mu_m : X^\pr_m \lra
X_m $ and $\mu_{\infty} :
X^\pr_{\infty}
\lra X_\infty$.

Fix next a multi-index ${\nu }= (\nu_i )
\in \NN^t$. We consider   the ``multi-contact" locus
\[\Cont^{\nu}(E)  \ = \ \big \{\, \gamma^\pr
\in
X^\pr_\infty \mid \ord_{\gamma^\pr}(E_i) = \nu_i \
\text{for} \ 1 \le i \le t \, \big \}. \]
These are again locally closed cylinders arising from
the corresponding subsets $\Cont^\nu(E)_m$ of
$X^\pr_m$.  The philosophy is that these loci can be
understood very concretely thanks to the fact that
$E$ is a simple normal crossing divisor: for example, 
$\Cont^{\nu}(E)_m$ (if non-empty) is a smooth
irreducible variety of codimension
$\sum \nu_i$ in $X^\pr_m$. 

Our first main result describes
the contact loci of $Y$ in terms of the
multi-contact loci associated to $E$.
In particular, we see that any 
irreducible component of the contact locus
$\Cont^p(\fra)$ associated to the ideal sheaf
$\fra$ on $X$ arises as the image of a unique such
multi-contact locus on
$X^\pr$:
\begin{theoremalpha} \label{Main.Thm.Ref} 
For every positive integer $p$,
we have a decomposition as a finite disjoint union
$$\Cont^p(\fra)=\bigsqcup_{\nu}\mu_{\infty}(\Cont^{\nu}(E)),$$
over those $\nu\in\NN^t$ such that $\sum_i\nu_ir_i=p$.
Each $\mu_{\infty}(\Cont^{\nu}(E))$ is a constructible cylinder
of codimension $\sum_i\nu_i(k_i+1)$. In particular, for every 
irreducible component $W$ of $\Cont^p(\fra)$
there is a unique
multi-index
$\nu$ as above such that $\Cont^{\nu}(E)$
dominates
$W$.  
\end{theoremalpha}

A related result appears as
Theorem~2.4 in \cite{denef0}.
Theorem~\ref{Main.Thm.Ref}  yields a quick proof of results of
\cite{Mu2} relating the dimensions of the arc spaces
$Y_\ell  \subseteq \ X_\ell$
 of $Y$ to the singularities of the
pair $(X, Y)$.   Keeping notation as above,
recall that the log-canonical threshold of
$\fra$ -- or of the pair
$(X,Y)$ -- is the rational number
\[  \lct(\fra) \ = \ \lct(X,Y) \ =_{\text{def}}
\min_{i} \Big \{ \frac{k_i + 1}{r_i} \Big \}. \] This
is an important invariant of the singularities of the
functions $f \in  \fra$, with smaller values of
$\lct(\fra)$ reflecting nastier singularities (see
\cite{DK}, \cite{ELSV}, \cite{kollar}).  Now consider an irreducible
component $V_\ell$ of $Y_\ell$. Then the
inverse image of $V_\ell$ in $X_{\infty}$ is contained in the
closure of an irreducible component $W$ of the contact
locus
$\Cont^p(Y)$ for some $p \ge \ell + 1$. Write $\nu =
(\nu_i)$ for the multi-index describing the cylinder
$\Cont^\nu(E)$ dominating $W$. If
$c =
\lct(X,Y)$ then
$k_i + 1
\ge c r_i$ for each of the  divisors $E_i$, and so one
finds from Theorem~\ref{Main.Thm.Ref}:
\begin{equation} \label{LCT.Bound.Eqn}
\begin{aligned} 
\codim\big( V_\ell, X_\ell\big) \ &\ge \ \codim\big(
\mu_{\infty}(\Cont^{\nu}(E)), X_{\infty}\big) \\
&= \ \sum \nu_i(k_i +1) \\ 
&\ge\ c \cdot \sum \nu_i r_i \\&= \ c\cdot p \\  &\ge \ c \big(
\ell + 1 \big).
\end{aligned}
\end{equation}
The reverse inequality being elementary for suitable
values of $\ell$, we deduce
\begin{corollaryalpha}[Musta\c{t}\v{a}, \cite{Mu2}]
\label{Cor.B.Ref} One
has
\[ \lct(X,Y) \ = \ \min_{\ell} \Big\{ \frac{\codim
(Y_\ell, X_\ell)}{\ell + 1} \Big\}. \]
\end{corollaryalpha} 
\noi A closer look at this argument identifies
geometrically the components of 
$Y_\ell$ having maximal possible dimension: they are the closures of the images
of multi-contact loci on $X^\pr$ involving the divisors
$E_i$ computing $\lct(\fra)$ (Corollary
\ref{Description.Max.Components}). We also recover the
main result of \cite{Mu1} relating irreducibility of
$Y_\ell$ to the singularities of $Y$ when $Y$ is a
local complete intersection.

Our second main result concerns the valuations defined
by irreducible cylinders $C \subset X_\infty$. Assuming
that $C$ does not dominate $X$, it determines a
valuation $\val_C$ on the function field $\CC(X)$ of $X$
by the rule:
\[   \val_C(f) \ = \ \ord_\gamma(f) \ \ \text{for
general
$\gamma \in C$}. \]
We prove  that any such valuation
 comes from a contact locus: 

\begin{theoremalpha} \label{Intro.Thm.C} Let $C
\subseteq X_{\infty}$ be an irreducible 
cylinder which does not dominate $X$, and denote by
$\val_C$ the corresponding valuation on $\CC(X)$. Then
there exist an ideal $\fra
\subseteq \OO_X$, an  integer $p > 0$ and an
irreducible component $W$ of $\Cont^p(\fra)$ such that
$\val_C=\val_W$.
Moreover ${\rm val}_C$ is a 
divisorial valuation, i.e. it agrees up to a constant
with the valuation given by the order of vanishing
along a divisor in a suitable blow-up of $X$.
\end{theoremalpha}
\noi The divisor in question is the exceptional divisor
in a weighted blow-up of $X^\pr$ along the intersection
of the relevant $E_i$ from Theorem
\ref{Main.Thm.Ref}. 

Theorem \ref{Main.Thm.Ref} is a consequence of a
result of Denef and Loeser describing the fibres of the
map $\mu_m : X_m^\pr \lra X_m$. This is also one of the
principal inputs to the theory of motivic integration,
so the techniques used here are certainly
not  disjoint from the methods of the third author in
\cite{Mu1} and \cite{Mu2}. However the present approach 
clarifies the geometric underpinnings of the results in
question, and it bypasses the combinatorial
complexities involved in manipulating motivic
integrals. 

In order to streamline the presentation, we collect
in Section 1 some basic results on cylinders in
arc spaces of smooth varieties.  
The proof of the above theorems are given in \S 2. Section 3
is devoted to some variants and further applications:
in particular, we complete the proof of Corollary
\ref{Cor.B.Ref} and explain how to recover some of the
results of \cite{Mu1}. We discuss also more general pairs
of the form $(X,\alpha\cdot Y-\beta\cdot Z)$, and interpret
the corresponding generalized log canonical thresholds in terms
of arc spaces. In particular, this yields a jet-theoretic description
of the multiplier ideals. As an example, 
we treat the case of monomial subschemes
in an affine space.

\bigskip

{\bf Acknowledgements}.
We are grateful to Fran\c{c}ois Loeser and Mihnea Popa for
useful discussions, and to the referee for his comments on the paper.
Research of the first author was partially supported by
NSF Grant 02-00278. Research of the second author was partially
supported by NSF Grant 01-39713. The third author served as a
Clay Mathematics Institute Long-Term Prize Fellow while this
research has been done.

\section{Cylinders in arc spaces}

We collect in this section some basic facts on cylinders.
Let $X$ be a fixed smooth, connected $d$-dimensional
complex variety. We have truncation morphisms
$\psi_m : X_{\infty}\longrightarrow X_m$, and 
$\tau_{m+1,m} : X_{m+1}\longrightarrow X_m$, such that
$\tau_{m+1,m}$ is locally trivial with fiber $\CC^d$. Note that
$X_{\infty}$ is the set of $\CC$-valued points of a scheme over $\CC$,
and we consider on $X_{\infty}$ the restriction of the Zariski topology
on this scheme. It is clear that 
this is equal with the inverse limit topology induced by the 
projections $\{\psi_m\}_m$.
Since every $\tau_{m+1,m}$ is flat, hence open, it follows that
$\psi_m$ is open for every $m$.

Recall that
a cylinder $C$ in $X_{\infty}$ is a subset of the form $\psi_m^{-1}(S)$,
for some $m$, and some constructible subset $S\subseteq X_m$.
We stress that all the points we consider in $X_{\infty}$
and $X_m$ are $\CC$-valued points.
The cylinders form an algebra of sets. It is clear that $C=\psi_m^{-1}(S)$
is open (closed, locally closed) if and only if so is $S$.
Moreover, we have $\overline{C}=\psi_m^{-1}(\overline{S})$,
so $\overline{C}$ is a cylinder.
Since the projection maps are locally trivial, it follows that
$C$ is irreducible if and only if $S$ is. In particular,
the decomposition of $S$ in irreducible components induces
a similar decomposition for $C$.
For every cylinder $C=\psi_m^{-1}(S)$, we put 
$\codim(C):=\codim(S,X_m)=(m+1)d-\dim(S)$.
If $C_1\subseteq C_2$ are cylinders, then 
$\codim(C_1)\geq \codim(C_2)$. Moreover, if $C_1$ and $C_2$
are closed and irreducible, then we have equality of codimensions 
if and only if $C_1=C_2$.
 
We say that a subset $T$ of $X_{\infty}$ is \emph{thin} if
there is a proper closed subscheme $Z$ of $X$ such that 
$T\subseteq Z_{\infty}$.

\begin{proposition}\label{not_thin}
If $C$ is a non-empty cylinder, then $C$ is not thin.
\end{proposition}

\begin{proof}
Suppose that $Z$ is a proper closed subset such that 
$C\subseteq Z_{\infty}$. In particular, we have
$\codim(C)\geq \codim(\psi_m^{-1}(Z_m))$ for every $m$.
On the other hand, it can be shown that 
$\lim_{m\to\infty}\codim(\psi_m^{-1}(Z_m))=\infty$
(see, for example, Lemma~3.7  in \cite{Mu1}).
This gives a contradiction.
\end{proof}

\begin{lemma}\label{infinite_intersection} 
If
$C_1\supseteq C_2\supseteq C_3\supseteq\ldots$ are non-empty cylinders
in $X_{\infty}$, then
$\bigcap_mC_m\neq\emptyset$.
\end{lemma}

\begin{proof}
We give a proof following \cite{batyrev}.
Note first that a similar assertion holds for a non-increasing
sequence of constructible subsets of a given variety, as we work
over an uncountable field. Moreover, it follows from definition
and Chevalley's theorem that for every cylinder $C\subseteq X_{\infty}$,
the image $\psi_p(C)$ is constructible for all $p$.

Consider the constructible subsets of $X$
$$\psi_0(C_1)\ \supseteq\
\psi_0(C_2) \ \supseteq\ \psi_0(C_3)\ \supseteq\ \ldots.$$ The above
remark shows that we can find $x_0\in\bigcap_m\psi_0(C_m)$. Replace now
$C_m$ by
$C'_m:=C_m\cap\psi_0^{-1}(x_0)$, which again form a non-increasing
sequence of non-empty cylinders.  Consider 
$$\psi_1(C'_1)\ \supseteq \ \psi_1(C'_2)\
\supseteq \ \psi_1(C'_3)\ \supseteq\ \ldots,$$ so that we can find
$x_1\in\bigcap_m\psi_m(C'_m)$. We replace now
$C'_m$ by $C''_m:=C_m\cap\psi_1^{-1}(x_1)$ and continue this way, to
get a sequence $\{x_p\}_p$, such that $\tau_{p+1,p}(x_{p+1})=x_p$
for every $p$, and such that $x_p\in\psi_p(C_m)$ for every $p$ and $m$.
This defines $x=(x_p)_p\in X_{\infty}$, and since each $C_m$ is a cylinder,
we see that $x\in C_m$ for every $m$.
\end{proof}

We consider now a proper, birational morphism 
$\mu : X'\longrightarrow X$ between smooth varieties, with
$\mu_m : X'_m\longrightarrow X_m$ the induced map on arc spaces.
The next lemma shows that $\mu_m$ is surjective for all $m\in\NN$.
In fact, we will show later the stronger fact that $\mu_{\infty}$ is 
surjective.

\begin{lemma}\label{finite_surj}
If $\mu$ is as above, then $\mu_m : X'_m\longrightarrow
X_m$ is surjective for every $m\in\NN$.
\end{lemma}

\begin{proof}
Let $Z\subset X$ be a proper closed subset such that $\mu$
is an isomorphism over $X\setminus Z$. It follows from the 
Valuative Criterion for properness that $X_{\infty}\setminus Z_{\infty}$
is contained in the image of $\mu_{\infty}$. On the other hand, for
every $u\in X_m$, we have $\psi_m^{-1}(u)\not\subseteq Z_{\infty}$
by Proposition~\ref{not_thin},
so $u\in {\rm Im}(\mu_m)$.
\end{proof}

We give now a criterion for the image of a cylinder by $\mu_{\infty}$
to be a cylinder. If $\mu$ is as above, we denote by $\psi_m$ and $\psi'_m$
the projections corresponding to $X$ and $X'$, respectively.

\begin{proposition}\label{image1}
With the above notation, if
$C'=(\psi'_m)^{-1}(S')$ is a cylinder so that $S'$ is a union of fibers of 
$\mu_m$, then $C:=\mu_{\infty}(C')$ is a cylinder in $X_{\infty}$.
\end{proposition}

\begin{proof}
Since $S:=\mu_m(S')$ is constructible, it is enough to show that
$C=\psi_m^{-1}(S)$. We prove that $C\supseteq\psi_m^{-1}(S)$,
the reverse inclusion being obvious. If $\gamma\in\psi_m^{-1}(S)$,
consider for every $p\geq m$, the cylinder 
$$D_p=(\psi'_p)^{-1}(\mu_p^{-1}(\psi_p(\gamma)))\subseteq X'_{\infty}.$$
It follows from Lemma~\ref{finite_surj} that $D_p\neq\emptyset$ for every $p\geq m$.
On the other hand, it is clear that $D_p\supseteq D_{p+1}$ for $p\geq m$,
so Lemma~\ref{infinite_intersection} 
shows that there is $\gamma'\in\bigcap_{p\geq m}D_p$.
Note that $\mu_{\infty}(\gamma')=\gamma$, while we have $\gamma'\in C'$,
as $S'$ is a union of fibers of $\mu_m$.
\end{proof}

We use this to strengthen the assertion in Lemma~\ref{finite_surj}.

\begin{corollary}\label{infinite_surj}
If $\mu : X'\longrightarrow X$
is a proper, birational morphism between smooth varieties, then
$\mu_{\infty}$ is surjective.
\end{corollary}

\begin{proof}
$C'=X'_{\infty}$ certainly satisfies the hypothesis of Lemma~\ref{image1},
so $C=\mu_{\infty}(X'_{\infty})$ is a cylinder. On the other hand,
we have seen in the proof of Lemma~\ref{finite_surj} that there is a proper
closed subset $Z\subset X$ so that $X_{\infty}\setminus Z_{\infty}
\subseteq C$. Therefore $X_{\infty}\setminus C$ is contained
in $Z_{\infty}$, so it is empty by Proposition~\ref{not_thin}.
\end{proof}

\smallskip

We recall now a theorem of Denef and Loeser which
will play a pivotal role in our arguments. 
Suppose that $\mu : X^\pr \lra X$ is a proper,
birational morphism of smooth varieties. As in the Introduction, we
introduce the relative canonical divisor
\[ K_{X^\pr / X} \ =_{\text{def}}\ \{\, \det (d\mu ) \, =
\, 0 \,\},\]
an effective Cartier divisor on $X'$, supported on the exceptional
locus of $\mu$. 
\begin{theoremnn}[Denef and Loeser, \cite{denef}]
Given an integer $e \ge 0$, consider the contact locus
\[ \Cont^e(K_{X^\pr/X})_m \ = \ \big\{\, \gamma^\pr \in
X^\pr_m  \mid \ord_{\gamma^\pr}(K_{X^\pr/X}) = e \,
\big \}. \]
If $m \ge 2e$ then $\Cont^e(K_{X^\pr/X})_m $ is
a union of fibres of $\mu_m : X_m^\pr \lra X_m$, each
of which is isomorphic to an affine space $\AAA^e$.
Moreover if \[
\gamma^\pr \ , \ \gamma^{\pr \pr} \ \in \ 
\Cont^e(K_{X^\pr/X})_m \]
 lie in the same fibre of
$\mu_m$, then they have the same image in
$X^\pr_{m-e}$. 
\end{theoremnn}
\noi In fact,  
Denef and Loeser show that  $\mu_m$ is a Zariski-locally trivial
$\AAA^e$-bundle over the image of 
$\Cont^e(K_{X^\pr/X})_m $  in $X_m$. 

\begin{remark} The most important point for our
purposes is the statement that the contact locus
$\Cont^e(K_{X^\pr/X})_m \subseteq	 X^\pr_m$ is a union
of fibres of $\mu_m$, and that all of  these fibres are
irreducible of dimension $e$. When
$\mu : X^\pr \lra X$ is the blowing-up of $X$ along a
smooth center this is readily checked by an explicit
calculation in local coordinates. This case  in turn
 implies the statement when $\mu$ is obtained as a
composition of such blow-ups. There would be no
essential loss in generality in limiting ourselves in
what follows to such a composition of nice blow-ups. 
\qed
\end{remark}

\begin{corollary}\label{image2}
Let $\mu : X'\longrightarrow X$
be a proper, birational morphism between smooth varieties.
If $C'\subseteq X'_{\infty}$ is a cylinder, then the closure
$\overline{\mu_{\infty}(C')}$ of its image is a cylinder.
Moreover, if there is $e\in\NN$ such that
$C'\subseteq
{\rm Cont}^e(K_{X'/X})$, then 
$\mu_{\infty}(C')$ is a cylinder. 
\end{corollary}

\begin{proof}
We use the above theorem of Denef and Loeser. 
Note first that in order to prove the second assertion,
we may assume that $C'$ is irreducible.
Let $e$ be the smallest $p$ such that $C'\cap {\rm Cont}^p(K_{X'/X})\neq
\emptyset$ (we have $e<\infty$ by Proposition~\ref{not_thin}). 
As $C'_{\circ}:=C'\cap {\rm Cont}^e
(K_{X'/X})$ is open and dense in $C'$, it follows that $\overline{\mu_{\infty}(C')}
=\overline{\mu_{\infty}(C'_{\circ})}$, so that it is enough to prove
the second assertion in the theorem.

Let $p$ and $T\subseteq X'_p$ be such that $C'=(\psi'_p)^{-1}(T)$,
and fix $m$ with $m\geq\max\{2e,e+p\}$. If $S\subseteq X'_m$ is the inverse
image of $T$ by the canonical projection to $X'_p$, Proposition~\ref{image1}
shows that it is enough to check that $S$ is a union of fibers of $\mu_m$.
This follows from Denef and Loeser's theorem: if 
$\delta_1\in S$ and $\delta_2\in X'_m$ are such that $\mu_m(\delta_1)
=\mu_m(\delta_2)$, then the first assertion in the theorem implies
$\delta_2\in {\rm Cont}^e(K_{X'/X})_m$, and the last assertion in the theorem
shows that $\delta_1$ and $\delta_2$ lie over the same element in $X'_p$,
hence $\delta_2\in S$.
\end{proof}

We show now that our notion of codimension 
for cylinders agrees with the usual one, defined in terms of the
Zariski topology on $X_{\infty}$. 

\begin{lemma}\label{codim1}
If $C\subseteq X_{\infty}$ is an irreducible cylinder,
and if $W\supseteq C$ is an irreducible closed subset of $X_{\infty}$, then
$W$ is a clinder.
\end{lemma}

\begin{proof}
It follows from the definition of the Zariski topology that
$W=\bigcap_{m\in\NN}W^{(m)}$, where $W^{(m)}=\psi_m^{-1}(\overline{\psi_m(W)})$.
Note that each $W^{(m)}$ is a closed irreducible cylinder, and we have
$C\subseteq W^{(m+1)}\subseteq W^{(m)}$ for every $m$. Therefore
$\codim\,W^{(m)}\leq\codim(C)$ for every $m$, and we deduce that 
$\codim\,W^{(m)}$ is eventually constant. This shows that there is
$m_0$ such that $W^{(m)}=W^{(m_0)}$ for all $m\geq m_0$;
hence $W=W^{(m_0)}$ is a cylinder.
\end{proof}

\begin{corollary}\label{codim2}
If $C\subseteq X_{\infty}$ is a cylinder, then
$\codim(C)$ is the codimension of $C$, as defined using the Zariski topology
on $X_{\infty}$. 
\end{corollary}

\begin{proof}
We may clearly assume that $C$ is closed and
irreducible. In this case, the above lemma
shows that every chain of closed irreducible subsets containing $C$
consists of cylinders, so the assertion is obvious.
\end{proof}

\begin{proposition}\label{disjoint_union}
Let $C\subseteq X_{\infty}$ be a cylinder.
If we have a countable disjoint union of cylinders
$\bigsqcup_{p\in\NN}D_p\subseteq C$,
whose complement in $C$ is thin, then  
$$\codim(C)=\min_{p\in\NN}\codim(D_p).$$
Moreover, if each $D_p$ is irreducible (or empty), then for every
irreducible component $W$ of $C$ there is a unique $p\in\NN$ such that
$D_p\subseteq W$, and $D_p$ is dense in $W$.
\end{proposition}

\begin{proof}
Let $Z\subset X$ be a proper closed subset
such that $C\subseteq Z_{\infty}\cup\bigcup_pD_p$.
We will use the fact that $\lim_{m\to\infty}\codim(\psi_m^{-1}(Z_m))=\infty$
(see Lemma~3.7 in \cite{Mu1}).

It is clear that we have $\codim(C)\leq\min_p\codim(D_p)$.
For the reverse inequality, choose $m$ such that $\codim(\psi_m^{-1}(Z_m))
>\min_p\codim(D_p)$. It follows from Lemma~\ref{infinite_intersection}
 that there is $r$
such that 
$$C\subseteq\bigcup_{p\leq r}D_p\cup\psi_m^{-1}(Z_m).$$
This clearly gives $\codim(C)\geq\min_p\codim(D_p)$.

Suppose now that every $D_p$ is irreducible, and let $W$ be an
irreducible component of $C$. Choose $m$ such that
$\codim(\psi_m^{-1}(Z_m))>\codim(W)$, and let $r$ be such that
$$C\subseteq\bigcup_{p\leq r}D_p\cup\psi_m^{-1}(Z_m).$$
It follows that there is $p\leq r$ such that $W\subseteq\overline{D_p}$.
Since $D_p$ is irreducible, we see that we also have $D_p\subseteq W$.
\end{proof}

\begin{remark}\label{motivic_measure}
Under the assumption of Proposition~\ref{disjoint_union}
one knows that
$$\mu_{\rm Mot}(C)=\sum_p\mu_{\rm Mot}(D_p),$$
where $\mu_{\rm Mot}$ is the motivic measure from \cite{denef}.
By taking the Hodge realization of the motivic measure,
the statement in the above proposition follows immediately.
However, we preferred to avoid this formalism.
\end{remark}

\section{Contact loci and valuations}

This section is devoted to the proof of our main
results. 
Keeping the notation established in the
Introduction, we give first
the statement and the proof of 
a more general version of Theorem \ref{Main.Thm.Ref}.
Consider a smooth variety
$X$ of dimension $d$ and a non-zero ideal sheaf $\fra
\subseteq \OO_X$ on $X$ defining a subscheme $Y
\subseteq X$. Fix a log resolution $\mu : X^\pr \lra X$
 of $\fra$, with $E = \sum_{i=1}^t E_i$ a simple normal
crossing divisor on $X^\pr$ such that
\begin{equation} \fra^\pr  \ \ = \ \fra \cdot
\OO_{X^\pr} \ = \ \OO_{X^\pr} \big(\, - \sum_{i=1}^t r_i
E_i \, \big) \quad \ ,
\ \quad K_{X^\pr / X} \ = \ \sum_{i=1}^t \, k_i E_i 
\end{equation}
for some integers $r_i, k_i \ge 0$. 
Note that this time we do not assume $r_i\geq 1$ if $k_i\geq 1$.
Given a multi-index
$\nu = (\nu_i) \in \NN^t$ we define the support of
$\nu$ to be
\[ \supp(\nu) \ = \ \big \{\,  i \in [1, t ] \, \mid \, 
\nu_i \ne 0 \, \big \}, \]
and we put 
\[ E_\nu \ = \ \bigcap_{i \in \supp(\nu)} \, E_i. \]
Thus $E_\nu$ is either empty or a smooth subvariety of
codimension $|\supp(\nu)|$ in $X'$. Without loss of
generality we will assume in addition that $E_\nu$ is
connected --- and hence irreducible --- whenever it is
non-empty.\footnote{We can always arrive at this
situation by a sequence of  blow-ups along
smooth centers. In fact, starting with an arbitrary log
resolution, first blow up the $d$-fold intersections of
of the $E_i$; then blow up
the $(d-1)$-fold intersections of their proper
transforms; and so on. At the end of this
process we arrive at a log resolution where the stated
condition is satisfied.}

Given a multi-index $\nu \in \NN^t$ and an integer $m
\ge \max_i{\nu_i}$,  we consider as in the Introduction
the ``multi-contact" loci
\[ \Cont^{\nu}(E)_m \ = \ \big \{\, \gamma^\pr
\in
X^\pr_m \mid \ord_{\gamma^\pr}(E_i) = \nu_i \
\text{for} \ 1 \le i \le t \, \big \},
\]
and the corresponding sets $\Cont^\nu(E) \subseteq
X^\pr_\infty$. Provided that $E_\nu \ne \emptyset$, 
a computation in local coordinates shows that
$\Cont^\nu(E)_m$ is a smooth irreducible locally closed
subset  of codimension $\sum \nu_i$ in $X^\pr_m$. (If
$E_\nu$ is empty, then so is $\Cont^\nu(E)_m$.) 
Recall that a subset of $X_{\infty}$ is called {\emph{thin}}
if it is contained in $Z_{\infty}$ for some proper closed subscheme
$Z\subset X$.

\begin{theorem}\label{main1}
For every positive integer $p$, we have a disjoint union
$$\bigsqcup_{\nu}\mu_{\infty}(\Cont^{\nu}(E))\subseteq \Cont^p(Y),$$
where the union is over those $\nu\in\NN^t$ such that $\sum_i\nu_ir_i=p$.
For every $\nu$ as above such that $\Cont^{\nu}(E)\neq\emptyset$, its
image $\mu_{\infty}(\Cont^{\nu}(E))$ is an irreducible cylinder
of codimension $\sum_i\nu_i(k_i+1)$. Moreover, the complement in $\Cont^p(Y)$ of
the above union is thin.
\end{theorem}

\begin{proof}
We use the Theorem of Denef and Loeser which
was stated in the previous section.
It is clear that we have a disjoint union
$$\bigsqcup_{\nu}\Cont^{\nu}(E)\subseteq\mu_{\infty}^{-1}(\Cont^p(Y)),$$
where $\nu$ varies over the set in the statement, and the complement
is contained in the union of those $(E_i)_{\infty}$ such that $r_i=0$.
For every $\nu$ we put $e=\sum_i\nu_ik_i$, so that
$\Cont^{\nu}(E)\subseteq\Cont^e(K_{X'/X})$.
Corollary~\ref{image2} implies that 
$\mu_{\infty}(\Cont^{\nu}(E))$ is a cylinder.

We show now that $\mu_{\infty}(\Cont^{\nu}(E))$
and $\mu_{\infty}(\Cont^{\nu'}(E))$ are disjoint if $\nu\neq\nu'$.
Indeed, otherwise there are $\gamma\in\Cont^{\nu}(E)$ and 
$\gamma'\in\Cont^{\nu'}(E)$ with $\mu_{\infty}(\gamma)=\mu_{\infty}(\gamma')$.
If $e$ and $e'$ correspond to $\nu$ and $\nu'$, respectively,
fix $m\geq\max\{2e,e+\nu_i,e+\nu'_i\}_i$. 
The first assertion in the Theorem of Denef and Loeser
gives $e=e'$, and the last assertion implies that $\gamma$ and $\gamma'$
have the same image in $X'_{m-e}$, a contradiction.

Since $\mu_{\infty}$ is surjective by Corollary~\ref{infinite_surj},
we get a disjoint union as in the theorem. 
Moreover, if $\nu$ is such that $E_{\nu}\neq\emptyset$ 
and if $m\gg 0$, then the Theorem of Denef and Loeser shows that the projection
$$\Cont^{\nu}(E)_m\longrightarrow\mu_m(\Cont^{\nu}(E)_m)$$
has irreducible, $e$-dimensional fibers (in fact, it is
locally trivial with fiber $\AAA^e$). Hence 
$$\dim\mu_m(\Cont^{\nu}(E)_m)=\dim\Cont^{\nu}(E)_m-e=(m+1)d-\sum_i\nu_i(k_i+1),$$
which completes the proof of the theorem.
\end{proof}

\begin{remark}
Note that if $\mu$ is as in the Introduction, i.e. if it is
an isomorphism over $X\setminus Y$, then the disjoint union
in the above theorem is finite. Moreover, it follows from the above proof
that in this case the union is equal with $\Cont^p(Y)$, and we recover also the
first statement of Theorem~\ref{Main.Thm.Ref}.
\end{remark}

\begin{remark}
As pointed out by the referee, the fact that 
$$\mu_{\infty}(\Cont^{\nu}(E))
\cap\mu_{\infty}(\Cont^{\nu'}(E))=\emptyset$$ for $\nu\neq\nu'$
can be seen also in an elementary way as follows. If $Z=\mu(E)$, then $\mu$
induces an isomorphism $X'\setminus E\longrightarrow X\setminus Z$.
By the Valuative Criterion for Properness, $\mu_{\infty}$
induces a bijection $X'_{\infty}\setminus E_{\infty}
\longrightarrow X_{\infty}\setminus Z_{\infty}$. As $\nu_i$, 
$\nu'_i\neq\infty$
for all $i$, we see that $\Cont^{\nu}(E)$, $\Cont^{\nu'}(E)
\subseteq X'_{\infty}\setminus E_{\infty}$, and we deduce our assertion.
\end{remark}

In the setting of Theorem~\ref{main1}, the disjoint decomposition
of $\Cont^p(Y)$ allows us to relate the irreducible components of
$\Cont^p(Y)$ with the multi-contact loci $\Cont^{\nu}(E)$. This is
a formal consequence of Proposition~\ref{disjoint_union}.

\begin{corollary}\label{compute_codim}
With the notation in Theorem~\ref{main1}, we have
$$\codim \, \Cont^p(Y)\ = \ \min_{\nu}\sum_i\nu_i(k_i+1),$$ 
where the
minimum is over all $\nu\in\NN^t$ such that $\sum_i\nu_ir_i=p$
and $E_{\nu}\neq\emptyset$. Moreover, for every irreducible
component $W$ of $\Cont^p(Y)$ there is a unique $\nu$ as above, such that
$W$ contains $\mu_{\infty}(\Cont^{\nu}(E))$ as a dense subset.
\end{corollary}

Note that this corollary allows us to describe the irreducible
components of each $\Cont^p(Y)$ which have minimal codimension
in terms of the numerical data
of the resolution. However, while Theorem~\ref{main1}
shows that all the irreducible
components are determined by the resolution, the ones
of codimension larger than $\codim(\Cont^p(Y))$
 seem to depend on
more than just the numerical data. In the next section
we will see that by considering a varying auxiliary scheme $Z\subset X$,
we can describe also other components of the contact loci.

\bigskip

We turn now to valuations of $\CC(X)$ associated to
cylinders in the arc space. Recall that if $C$ is an irreducible cylinder
in $X_{\infty}$, then we have defined a valuation ${\rm val}_C$
as follows. Note first that if $\gamma\in X_{\infty}$,
and if $f$ is a rational function on $X$ defined in a neighbourhood
of $\psi_0(\gamma)$, then ${\rm ord}_{\gamma}(f)$ is well defined.
If the domain of $f$ intersects $\psi_0(C)$, then 
${\rm val}_C(f):={\rm ord}_{\gamma}(f)$, for general $\gamma\in C$. 
Note that this is a non-negative integer by Proposition~\ref{not_thin}.
Since $C$ is irreducible, it follows that ${\rm val}_C(f)$
is well-defined and can be extended to a valuation of the function field of $X$.
We will assume from now on that $C$ does not dominate $X$, so
${\rm val}_C$ is nontrivial (and discrete).

It is clear that if $C_1\subseteq C_2$ are cylinders as above,
then ${\rm val}_{C_1}(f)\geq{\rm val}_{C_2}(f)$ for every rational
function $f$ whose domain intersects $\psi_0(C_1)$. If moreover,
$C_1$ is dense in $C_2$, then ${\rm val}_{C_1}={\rm val}_{C_2}$.
Suppose now that $\mu : X'\longrightarrow X$ is a proper, birational
morphism of smooth varieties, and that $C'\subseteq X'_{\infty}$
is an irreducible cylinder which does not dominate $X'$. It follows
from Corollary~\ref{image2} that $C:=\overline{\mu_{\infty}(C')}$
is a cylinder, and we clearly have
${\rm val}_{C'}={\rm val}_C$.

\begin{example}\label{cylinder_of_divisor}
A divisorial valuation of $\CC(X)$ is a discrete valuation 
associated to a prime divisor on some normal variety $X'$
which is birational to $X$. If the valuation has center on $X$
(see below), then
there are infinitely many cylinders
in $X_{\infty}$ whose corresponding valuations agree up to
a constant with the given valuation. To see this
let $D$ be a prime divisor on $X'$, and let ${\rm val}_D$
be the corresponding valuation. Saying that the valuation
has center on $X$ means that
we may assume 
that we have a proper birational morphism $\mu : X'\longrightarrow X$,
that $X'$ is smooth, and that $D$ is smooth on $X'$. For $p\geq 1$,
let $C'_p=\Cont^p(D)$, so ${\rm val}_{C'_p}=p\cdot{\rm val}_D$.
If we put $C_p=\overline{\mu_{\infty}(C'_p)}$, then Corollary~\ref{image2}
shows that $C_p$ is an irreducible cylinder, and ${\rm val}_{C_p}=
p\cdot{\rm val}_D$.  
\end{example}

In the remainder of this section we show that conversely,
every valuation defined by a cylinder in $X_{\infty}$ is
(up to a constant multiple)
a divisorial valuation. We consider first the case when the
cylinder is an irreducible component of a contact locus.
In this case, the assertion is a corollary of Theorem~\ref{main1}.

\begin{corollary}\label{irred_comp}
Let $Y\subset X$ be proper closed subscheme, and let 
$W$ be an irreducible component of $\Cont^p(Y)$ for some $p\geq 1$.
Let $\nu\in\NN^t$ be the multi-index given by Corollary~\ref{compute_codim},
so that $\mu_{\infty}(\Cont^{\nu}(E))$ is contained and dense in $W$.
If $D$ is the exceptional divisor of the weighted blowing-up of
$(X',E)$ with weight $\nu$, then ${\rm val}_W=q\cdot{\rm val}_D$,
where $q={\rm gcd}\{\nu_i\vert i\in {\rm supp}(\nu)\}$.
\end{corollary}

\begin{proof}
Recall the definition of the weighted blowing-up with weight $\nu$. 
Let $s={\rm lcm}\{\nu_i\vert i\in {\rm supp}(\nu)\}$.
If $T_{\nu}\subset X'$ is the closed subscheme defined by
$\sum_{i\in {\rm supp}(\nu)}\OO(-\frac{s}{\nu_i}\cdot E_i)$,
then the weighted blowing-up of $(X',E)$ with weight $\nu$
is the normalized blowing-up of $X'$ along $T_{\nu}$.
There is a unique prime divisor on this blowing-up which
dominates $T_{\nu}$; this is $D$.

It follows from our choice of $\nu$ that
${\rm val}_W={\rm val}_{\Cont^{\nu}(E)}$. Let $g : X''\longrightarrow X'$
be a proper birational map which factors through the above blowing-up
and which satisfies the requirements for a log resolution. Note that
we may consider $D$ as a divisor on $X''$. We apply 
Theorem~\ref{main1} for $g$. 

Let $C$ be the multi-contact locus
 of all arcs on $X''_{\infty}$ with order $q$
along $D$, and order zero along all the other 
divisors involved. Well-known results about weighted
blow-ups show that the coefficient of $D$ in $g^{-1}(E_i)$ is 
$\nu_i/q$, and the coefficient of $D$ in $K_{X''/X'}$
is $-1+\sum_i\nu_i/q$. We see that $g_{\infty}(C)\subseteq\Cont^{\nu}(E)$.
Moreover, both these cylinders are irreducible and have the
same codimension, as $\codim\,g_{\infty}(C)=q\cdot\sum_i\nu_i/q=
\codim(\Cont^{\nu}(E))$. This gives ${\rm val}_W={\rm val}_C=q\cdot{\rm val}_D$. 
\end{proof}

We show now that in fact, we can describe all
 valuations given by cylinders using contact loci.

\begin{theorem}\label{main2}
If $C$ is an irreducible cylinder in $X_{\infty}$ which does not
dominate $X$, then there is a proper closed subscheme $Y\subset X$,
a positive integer $p$, and an irreducible component $W$ of $\Cont^p(Y)$
such that ${\rm val}_C={\rm val}_W$.
In particular, ${\rm val}_C$ is equal, up to a constant, with a divisorial 
valuation.
\end{theorem}

\begin{proof}
We have to prove only the first assertion, the second one
follows from this and Corollary~\ref{irred_comp}.
By replacing $C$ with its closure, we may assume that $C$ is closed.
Moreover, it is enough to prove our assertion in the case when
$X=\Spec(A)$ is affine. Recall that a graded sequence of ideals
is a set of ideals $\fra_{\bullet}=\{\fra_p\}_{p\geq 1}$ such that
$\fra_p\cdot\fra_q\subseteq\fra_{p+q}$ for all $p$ and $q$
(see \cite{lazarsfeld} for more on this topic). Since ${\rm val}_C$
is a valuation which is non-negative on $A$, if we define
$$\fra_p:=\{f\in A\vert {\rm val}_C(f)\geq p\},$$
then $\fra_{\bullet}$ is a graded sequence of ideals.
Note that since $C$ does not dominate $X$, we have $\fra_p\neq (0)$ for every $p$.

Starting with a graded sequence of ideals $\fra_{\bullet}$ as above,
we get a sequence of closed cylinders as follows: for every $p\geq 1$, let
$$W_p=\{\gamma\in X_{\infty}\vert\ord_{\gamma}(f)\geq p\,\,{\rm for}\,{\rm every}\,
f\in\fra_p\}.$$
Since $\fra_{\bullet}$ is a graded sequence, we have $\fra_p^q\subseteq
\fra_{pq}$, so that $W_{pq}\subseteq W_p$ for every $p$, $q\geq 1$. 
Note that in our case, it follows from definition that $C\subseteq W_p$
for all $p$. Moreover, since $\fra_p\neq(0)$, we see that $W_p$
does not dominate $X$, for any $p$. We put $C_m:=W_{m!}$
so that $C\subseteq C_{m+1}\subseteq C_m$ for every $m\geq 1$.

We claim that we can choose irreducible components $C'_m$ of $C_m$
such that $C\subseteq C'_{m+1}\subseteq C'_m$ for all $m$. It is clear
that for every $m$ we can choose irreducible components $C_{i,m}$ of $C_i$
for $i\leq m$, so that
$$C\subseteq C_{m,m}\subseteq C_{m-1,m}\subseteq\ldots\subseteq C_{1,m}.$$
As every cylinder has finitely many irreducible components, there
is an irreducible component $C'_1$ of $C_1$ such that $C'_1=C_{1,m}$
for infinitely many $m$. Similarly, there is an irreducible component
$C'_2$ of $C_2$, such that $C'_2\subseteq C'_1$, and such that 
$C'_2=C_{2,m}$ for infinitely many $m$. Continuing in this way, we deduce our claim.

Note that we have $\codim(C'_1)\leq\codim(C'_2)\leq\ldots\leq\codim(C)$.
Therefore there is $q$ such that $\codim(C'_m)=\codim(C'_q)$ for every $m\geq q$,
hence $C'_m=C'_q$, as all $C'_i$ are irreducible closed cylinders.
Let $Y$ be the closed subscheme defined by $\fra_{q!}$, and
let $\tau:=\min\{\ord_{\gamma}(f)\vert\gamma\in C'_q, f\in\fra_{q!}\}$.
It is clear that $\tau\geq q!$ and that $C'_q$ is the closure
of an irreducible component of $\Cont^{\tau}(Y)$, so that in order
to finish the proof, it is enough to show that ${\rm val}_C={\rm val}_{C'_q}$.

Since $C\subseteq C'_q$, we have ${\rm val}_C\geq {\rm val}_{C'_q}$ on $A$.
Fix $f\in A$, and let us show that $m:={\rm val}_C(f)\leq {\rm val}_{C'_q}(f)$.
By multiplying $f$ with $g$ such that ${\rm val}_{C}(g)>0$, we may
assume $m\geq 1$. Moreover, by taking a suitable power of $f$,
we may assume that $m=p!$ for some $p\geq q$.
By definition we have $f\in\fra_m$,
hence 
$$C'_q=C'_p\subseteq W_{p!}\subseteq
\{\gamma\in X_{\infty}\vert\ord_{\gamma}(f)
\geq m\}.$$
This gives ${\rm val}_{C'_q}(f)\geq m$, and completes the proof of the theorem. 
\end{proof}

\begin{remark}
It follows from Theorem~\ref{main2} that to each irreducible cylinder $C$ 
which does not dominate $X$ we may
associate a (unique) divisor $D$ over $X$ such that ${\rm val}_C=
\lambda\cdot{\rm val}_D$ for some $\lambda>0$
(of course, we identify two divisors over $X$ which give 
the same valuation). This map is obviously
not injective, but it is surjective by Example~\ref{cylinder_of_divisor}.

Suppose now that $Y\subset X$ is a fixed proper closed subscheme.
It would be interesting to understand which divisors appear from
irreducible components of contact loci of $Y$. One can consider this
as an embedded version of Nash's problem (see \cite{nash}).
If $Y$ is a variety, and if $\psi_0 : Y_{\infty} \longrightarrow Y$ is
the canonical projection, then
Nash described an injective map from
the set of irreducible components of $\psi_0^{-1}(Y_{\rm sing})$
to the set of ``essential'' exceptional divisors over $Y$ (divisors which
appear in every resolution of $Y$). He conjectured that this 
map is surjective, but a counterexample has been recently found in 
\cite{ishiikollar}.

In the next section we will use Theorem~\ref{main1} to describe certain
divisors which are associated to distinguished irreducible components of
the contact loci of $Y$, namely the divisors which compute generalized
versions of the log canonical threshold. 
\end{remark}

\section{Applications to log canonical thresholds}

We discuss here some  applications of
Theorem~\ref{main1}. As before, $X$ is a smooth irreducible
variety of dimension $d$,  $Y \subseteq X$ is a
subscheme defined by a non-zero ideal sheaf $\fra
\subseteq
\OO_X$, and  $\mu : X^\pr \lra
X$ is a log resolution  of $(X,Y)$ with
\[ \fra \cdot \OO_{X^\pr} \ = \ \OO_{X^\pr}(-\sum  r_i
E_i) \ \ , \ \ K_{X^\pr / X} \ =  \ \sum k_i E_i. \]

\noi \textit{Arc spaces of subschemes.} We indicate how
to recover some of the results of \cite{Mu1} and
\cite{Mu2} relating singularities of the pair $(X,Y)$
to the properties of the arc spaces $Y_\ell \subseteq
X_\ell$. 

Note that for every irreducible cylinder $C\subseteq X_{\infty}$
there is a subcylinder $C_0\subseteq C$ which is open in $C$, such that
$\ord_{\gamma}(\fra)$ is constant for $\gamma\in C_0$.
We denote this positive integer by $\ord_C(\fra)$ or $\ord_C(Y)$.
It is clear that if $p=\ord_C(Y)$, then $C$ is contained in
the closure of an irreducible component of $\Cont^p(Y)$.

Suppose now that $V_\ell$ is an irreducible component of $Y_\ell$,
and let $V=\psi_{\ell}^{-1}(V_\ell)\subseteq X_{\infty}$.
 Note that $V$ is an irreducible component of
\begin{equation} \label{Inv.Im.Comp.Eqn}
 \psi_{\ell}^{-1} \big ( Y_\ell
\big)
\  = \ \Cont^{\ge (\ell + 1)} \big( Y \big) \ 
 =_{\text{def}} \ \big \{ \, \gamma \in X_{\infty} \mid \
\ord_{\gamma}(Y) \ge \ell + 1 \, \big \}.
  \end{equation} 
Let $p=\ord_V(Y)$, so $p\geq \ell+1$. Note that
we might have strict inequality (Example~\ref{Big.Order.Example}).
In any case, $V$ is the closure of an irreducible component
of $\Cont^p(Y)$. Conversely, the closure of every irreducible component
of $\Cont^p(Y)$ is the inverse image of an irreducible component of $Y_p$.
Therefore our analysis of the contact loci $\Cont^p(Y)$
gives complete control over the arc spaces of $Y$.

\begin{example}\label{Big.Order.Example} Let $X =
\AAA^1$ with coordinate
$t$, and let $Y\subset X$ be defined by $(t^e) \subseteq \CC[t]$
for a fixed integer $e \ge 2$.  
Then $V_1=Y_1$ is the irreducible subset of $X_1$
consisting of arcs centered at the origin.
However,
$\ord_V(Y)=e>1$.
\end{example}

We start by completing the proof of Corollary
\ref{Cor.B.Ref} from the Introduction, namely 
we prove

\begin{corollary}\label{lct}
The log canonical threshold of $(X,Y)$ is given by
\[ \lct(X,Y) \ = \ \min_{\ell} \Big\{ \frac{\codim
(Y_\ell, X_\ell)}{\ell + 1} \Big\}. \]
\end{corollary}

\begin{proof}
We saw in the Introduction that Theorem
\ref{Main.Thm.Ref}  implies the inequality
\begin{equation} \codim ( Y_\ell, X_\ell) \ \ge \ c(\ell
+ 1)\end{equation}
 where $c = \lct(X,Y)$, so it remains only to
prove the reverse inequality for suitable $\ell$.  But
this is immediate. In fact,   it follows from the
definition of $\lct(X,Y)$ that there exists an index
$i$   --- say $i = 1$ --- for which
$k_1 + 1 = c r_1$. Let $\nu = (1,0,0, \ldots,0)$ be the
multi-index with $\nu_1 = 1$ and $\nu_i = 0$ for $i >
1$. 
It follows from Theorem~\ref{main1} that $\mu_{\infty}(\Cont^{\nu}(E))$
is a subcylinder of $\Cont^{r_1}(Y)$ of codimension $k_1+1$.
If ${\ell}=r_1-1$, then the closure of this subcylinder can be written as
$\psi_{\ell}^{-1}(V)$ for some closed subset $V\subseteq Y_{\ell}$, with
$\codim(V,X_{\ell})=c(\ell+1)$. The first part of the proof implies that
$V$ must be an irreducible component of $Y_{\ell}$, so we are done.
\end{proof}

The argument just completed leads to an explicit
description of the components of $Y_\ell$ having
maximal possible dimension. Keeping notation as before,
let us say that one of the divisors
$E_i
\subseteq X^\pr$ \textit{computes the log canonical
threshold} of $(X,Y)$  if $\lct(X,Y) = \frac{k_i +
1}{r_i}$. Note that in  general there may be several
divisors $E_i$ that compute this threshold.

\begin{corollary} \label{Description.Max.Components}
Let $V_\ell \subseteq Y_\ell$ be an irreducible
component of maximal possible dimension, i.e. with
\[ \codim (V_\ell, X_\ell) \ = \ (\ell + 1)
\cdot \lct(X,Y)  ,
\] and let $V = \psi_\ell^{-1} (V_\ell)$ be the
corresponding subset of $X_\infty$. Then
$\ord_V(\fra) =
\ell + 1$ and
$V$ is dominated by a multi-contact locus $\Cont^\nu(E)$
where
\[ \nu_i \ne 0 \ \ \Longrightarrow   \ \ E_i \text{
computes the log canonical threshold of  $(X,Y)$.}\]
Conversely the image of any such multi-contact locus
$\Cont^\nu(E)$ determines a component of $Y_{\ell}$
of maximal possible dimension.
\end{corollary}

\begin{proof} 
Write $c = \lct(X,Y)$. We return to the proof of
Corollary \ref{Cor.B.Ref} from the Introduction.
Specifically,
$V$ is contained in the closure of some
irreducible component
$W$ of
$\Cont^p(\fra)$, which in turn is dominated by some
multi-contact locus
$\Cont^\nu(E)$.  But
$\codim(V_\ell, X_\ell) = c(\ell + 1)$ by hypothesis,
and hence
 equality must hold in all the inequalities
appearing in equation (\ref{LCT.Bound.Eqn}) from the
Introduction. Therefore 
\[ p \ = \ \sum \nu_i r_i  \ = \ (\ell + 1).
\] Moreover since in any event $k_i + 1 \ge c r_i$ for
all $i$, we deduce that $\nu_i (k_i + 1) = c \nu_i r_i$
for each index
$i$. In particular, if $\nu_i \ne 0$ then $E_i$
computes the log-canonical threshold of $(X,Y)$. We
leave the converse to the reader.
\end{proof}

Similar arguments allow one to eliminate motivic
integration from  the main results of the third author
in 
\cite{Mu1}.  For example:
\begin{corollary} \textnormal{\rm (Musta\c{t}\v{a},
\cite{Mu1}).} Let $Y \subseteq X$ be a
reduced and irreducible locally complete intersection
subvariety of codimension $f$. Then the arc-space
$Y_\ell$ is irreducible for all $\ell \ge 1$ if and only
if
$Y$ has at worst rational singularities. 
\end{corollary}

\begin{proof} Following \cite{Mu1} let $\mu: X^\pr \lra
X$ be a log resolution of $(X, Y)$ which dominates the
blowing-up of $X$ along $Y$. Write $E_1$ for the
(reduced and irreducible) exceptional divisor created by
this blow-up, so that
\[ r_1 \ = \ 1 \ \ \ , \ \ \ k_1  \ = \ f-1, \]
and otherwise keep notation and assumptions as above. 
It is established in \cite{Mu1}, Theorem 2.1 and Remark
2.2, that $Y$ has at worst rational singularities if
and only if $k_i \ge f r_i$ for every index $i \ge 2$.
So we are reduced to showing:
\begin{equation}
Y_\ell \text{ is irreducible for all $\ell \ge 1$} \quad
\Longleftrightarrow \quad k_i  \ \ge f  r_i\,\,{\rm for}\, i
\ge 2. \tag{*}
\end{equation}
(It is  the proof of this equivalence  in
\cite{Mu1} that uses motivic integration.)

Assuming that $k_i \ge f r_i$ for $i \ge 2$ we show
that each $Y_\ell$ is irreducible. Note to begin with
that $Y_\ell$ has one ``main component"
$Y_{\ell}^{\text{main}}$, namely the closure of the
arc-space
$(Y_{\text{reg}})_\ell$: this component is dominated by
the multi-contact locus
$\Cont^{(\ell+1, 0 , 0, \ldots, 0)}(E)$ described by
the multi-index $\nu = (\ell+1, 0 , 0 , \ldots, 0)$.
Suppose for a contradiction that there is a further
component $V_\ell$ of $Y_\ell$. In the usual way, $V
= \psi_{\ell}^{-1}(V_\ell)$ lies in  the closure of
an irreducible component $W$ of $\Cont^p(Y)$ for
 $p \ge \ell + 1$, which via Corollary~\ref{irred_comp}
is dominated by a multi-contact locus $\Cont^\nu(E)$
for some $\nu = (\nu_i) \ne (\ell+1, 0 ,\ldots, 0)$.
Since
$Y \subseteq X$ is a local complete intersection of
codimension $f$, we have in any event
$\codim (W, X_\ell)   \le   (\ell + 1)\cdot f$.
In view of the hypothesis $k_i \ge fr_i$ for $i \ge 2$
we then find the inequalities:
\begin{align*}
(\ell + 1) \cdot f \ &\ge \ \codim(W) \\
&= \ \sum_{i \ge 1} \nu_i ( k_i + 1) \\
&= \ \nu_1 \cdot f  \, + \, \sum_{i \ge 2} \nu_i (k_i
+ 1)
\\ &\ge
\ f
\cdot  \sum_{i \ge 1} \nu_i r_i \, + \, \sum_{i \ge 2}
\nu_i
\\ &= \ fp \, + \, \sum_{i \ge 2} \nu_i. 
\end{align*}
But since $p \ge \ell + 1$ this forces $\nu_i = 0$ for
$i \ge 2$, a contradiction. 

Conversely, suppose that $k_i < f r_i$
for some $i \ge 2$: say $k_2 \le f r_2 -1$. Setting
$\nu = (0,1,0  ,\ldots,0)$ and $\ell = r_2 -1$, the
multi-contact locus $\Cont^\nu(E)$ maps to an
irreducible set $W_\ell \subseteq Y_\ell $ with
\[ \codim(W_\ell, X_\ell)  \ \le  \ (\ell + 1)f \ = \
\codim( Y_{\ell}^{\text{main}}, X_\ell),\]
and therefore $Y_\ell$ cannot be irreducible. 
\end{proof}

\bigskip

\noi \textit{Generalized log canonical thresholds.}
We extend now the above results to take into account
also an extra scheme $Z$. We fix two proper closed
subschemes $Y$, $Z\subseteq X$
defined by the ideal sheaves $\fra$ and $\frb$, respectively. 
We use the previous notation
for $\mu : X'\longrightarrow X$, but this time we assume
that $\mu$ is a log resolution for $Y\cup Z$. We write
$\frb\cdot\OO_{X'}=\OO_{X'}(-\sum s_iE_i)$. Having fixed 
also $\beta\in\QQ_+$, we define the log canonical threshold
$\lct(X,Y;\beta\cdot Z)$ to be the largest $\alpha\in\QQ_+$
such that $(X,\alpha\cdot Y-\beta\cdot Z)$ is log canonical, i.e.
$k_i+1\geq\alpha\cdot r_i-\beta\cdot s_i$ for all $i$. Therefore
\[  \lct(X,Y;\beta\cdot Z) =
\min_{i} \Big \{ \frac{k_i + 1+\beta\cdot s_i}{r_i} \Big \}. \]
It is standard to see that the definition does not depend on
the particular log resolution (see \cite{kollar}). 

The following corollary is a generalization of Corollary~\ref{lct}
to this setting. This time we state the formula in terms of the contact
loci, and leave the corresponding statement in terms of arc spaces
to the reader.

\begin{corollary}\label{lct2}
We have 
$$\lct(X,Y;\beta\cdot Z)=\min_C\Big\{\frac{\codim(C)+\beta\cdot\ord_C(Z)}
{\ord_C(Y)}\Big\}$$
where $C$ runs over the irreducible cylinders in $X_{\infty}$
which do not dominate $X$; in fact, it is enough to let $C$ run only 
over the irreducible components of $\Cont^p(Y)$, for $p\geq 1$.
\end{corollary}

\begin{proof}
Let $c=\lct(X,Y;\beta\cdot Z)$.
We show first that if $C$ is an irreducible cylinder in $X_{\infty}$,
then $\codim(C)\geq c\cdot\ord_C(Y)-\beta\cdot\ord_C(Z)$. 
Let $p=\ord_C(Y)$, and let $C_0\subseteq C$ be an open subcylinder such that
$C_0\subseteq\Cont^p(Y)$. Let $W$ be an irreducible component of
$\Cont^p(Y)$ containing $C_0$, and let $\nu$ be the multi-index
corresponding to $W$ by Corollary~\ref{compute_codim}. We have
$$\codim(C)=\codim(C_0)\geq\codim(W)=\sum\nu_i(k_i+1)$$
$$\geq\sum_i\nu_i(c\cdot r_i-\beta\cdot s_i)=c\cdot p-\beta\cdot\ord_W(Z)
\geq c\cdot p-\beta\cdot\ord_C(Z),$$
as required. 

In order to finish the proof, it is enough
to show that there is $p\geq 1$ and an irreducible component
$V$ of $\Cont^p(Y)$ such that $\codim(V)=c\cdot p-\beta\cdot\ord_V(Z)$.
For this, let $i$ be such that $c=(k_i+1+\beta\cdot s_i)/r_i$,
and take $p=r_i$ and $\nu$ such that $\nu_i=1$ and $\nu_j=0$ if $j\neq i$.
We may take $V$ to be the closure of $\mu_{\infty}(\Cont^{\nu}(E))$
in $\Cont^p(Y)$.
\end{proof}

We have a similar generalization of
Corollary~\ref{Description.Max.Components}. Again, we phrase
the result in terms of irreducible components of contact loci of $Y$.
We say that a divisor $E_i$ computes $\lct(X,Y;\beta\cdot Z)$ if 
$\lct(X,Y;\beta\cdot Y)=(k_i+1+\beta\cdot s_i)/r_i$.

\begin{corollary}\label{Description.Max.Components2}
Let $W$ be an irreducible component of $\Cont^p(Y)$, with $p\geq 1$, such that
\begin{equation}\label{extremal}
\codim(W)=p\cdot\lct(X,Y;\beta\cdot Z)-\beta\cdot\ord_W(Z).
\end{equation}
Then $W$ is dominated by a multi-contact locus $\Cont^{\nu}(E)$ where
 \[ \nu_i \ne 0 \ \ \Longrightarrow   \ \ E_i \text{
computes $\lct(X,Y;\beta\cdot Z)$.}\]
Conversely, the image of any such multi-contact locus determines
an irreducible component of $\Cont^p(Y)$ as above.
\end{corollary}

\begin{proof}
The proof is similar to that of Corollary~\ref{Description.Max.Components},
so we omit it.
\end{proof}

\begin{remark}
Recall that we have defined in the previous section a map
from the irreducible components of $\Cont^p(Y)$ with $p\geq 1$
to the divisors over $X$. One can interpret the above corollary
as saying that for every $Z$ and $\beta$, every divisor $D$ which computes
$\lct(X,Y;\beta\cdot Z)$ is in the image of this map. Moreover, 
the irreducible components which correspond to these divisors are 
precisely the ones satisfying (\ref{extremal}). 
\end{remark}

\begin{remark}\label{multiplier_ideals}
Suppose that $X=\Spec(A)$ is affine, that $Y$ is defined by the ideal 
$\fra\subset A$, and that $Z$ is defined by a principal ideal $(f)$.
It follows from definition that $\lambda<\lct(X,Y;Z)$
if and only if $f\in\II(X,\lambda\cdot\fra)$, the multiplier ideal of
$\fra$ with coefficient $\lambda$ (we refer to \cite{lazarsfeld}
for the basics on multiplier ideals).
One can therefore interpret Corollary~\ref{lct2}
as giving an arc-theoretic interpretation of multiplier ideals.
\end{remark}

\bigskip

We end with an example: monomial ideals in the polynomial ring.
For a non-zero monomial ideal $\fra$ in the polynomial ring
$R=\CC[T_1,\ldots,T_d]$, the Newton polyhedron of $\fra$, denoted by
$P_{\fra}$, is the convex hull in $\RR^n$ of the set $\{{\bf
u}\in\NN^d\mid
 T^{\bf u}\in \fra\}$. For ${\bf u}=({\bf u}_i)\in\NN^d$, we use the notation
$T^{\bf u}=\prod_iT_i^{{\bf u}_i}$.

\begin{proposition}\label{monomial}
Let $X=\AAA^d$, and $Y$, $Z\hookrightarrow X$ subschemes defined by
nonzero monomial ideals $\fra$ and $\frb\subseteq R$, respectively, where
$R=\CC[T_1,\ldots,T_d]$.
For every $\alpha$, $\beta\in\QQ_+$, with $\alpha\neq 0$,
we have $(X, \alpha\cdot Y-\beta\cdot Z)$ log canonical
if and only if
$$\beta\cdot P_{\frb}+{\bf e}\subseteq \alpha\cdot P_{\fra},$$
where ${\bf e}=(1,\ldots,1)$. 
\end{proposition}

\begin{proof}
For every ${\bf q}\in\NN^d$, consider the multi-contact
locus
$C_{\bf q}\subseteq X_{\infty}$, consisting of those arcs
with order ${\bf q}_i$ along the divisor defined by $T_i$, for all $i$.
Every contact locus of $Y$ is a union of such cylinders.
It follows from Corollary~\ref{lct2} 
that
$(X, \alpha\cdot Y-\beta\cdot Z)$ is log canonical if and only if

\begin{equation}\label{condition_monomial}
\beta\cdot\ord_{C_{\bf q}}(Z)\geq \alpha\cdot \ord_{C_{\bf
{q}}}(Y) -\codim(C_{\bf q}),
\end{equation}
for every ${\bf {q}}$.

It is easy to see that $\codim(C_{\bf q})=\sum_i{\bf q}_i$.
Moreover, it is clear that we have $\ord_{C_{\bf q}}(Y)
=\inf\{\sum_i{\bf u}_i{\bf q}_i\vert{\bf u}\in P_{\fra}\}$, and a
similar formula for $\ord_{C_{\bf q}}(Z)$. If we consider the
dual polyhedron of $P_{\fra}$,
$$P_{\fra}^{\circ}:=\{{\bf p}\in\RR^d\vert\sum_i{\bf p}_i{\bf u}_i\geq 1\,
\text{for all}\, {\bf u}\in P_{\fra}\},$$
then we see that $\ord_{C_{\bf q}}(Y)\geq m$ if and only if $(1/m){\bf q}
\in P_{\fra}^{\circ}$.

Condition~(\ref{condition_monomial}) then says that for every
${\bf q}\in\NN^d$ and every $m\in\NN^*$,
 such that $(1/m){\bf q}\in P_{\fra}^{\circ}$,
we have
$$\beta/m\cdot\inf_{{\bf v}\in P_{\frb}}
\left(\sum_i{\bf v}_i{\bf q}_i\right)+\sum_i{\bf q}_i/m\geq\alpha.$$
Since $P_{\fra}^{\circ}$ is a rational polyhedron, this is the same as saying
that for every ${\bf q}'\in P_{\fra}^{\circ}$, we have
$$\beta\cdot\inf_{{\bf v}\in P_{\frb}}\left(\sum_i{\bf q}'_i{\bf v}_i\right)
+\sum_i{\bf q}'_i\geq \alpha.$$

Since $(P_{\fra}^{\circ})^{\circ}=P_{\fra}$ and $\alpha>0$,
the above condition is equivalent
to $\beta\cdot P_{\frb}+{\bf e}\subseteq \alpha\cdot P_{\fra}$.
\end{proof}

\begin{corollary}{\rm (Howald, \cite{howald})}
If $X=\AAA^d$ 
and $\fra\subseteq R=\CC[T_1,\ldots, T_d]$ is a nonzero monomial ideal,
then for every $\alpha\in\QQ_+^*$, the multiplier ideal of
$\fra$ with coefficient $\alpha$ is given by
$$\II(X, \alpha\cdot \fra)=(T^{\bf u}\mid {\bf u}+{\bf e}\in {\rm Int}
(\alpha\cdot P_{\fra})).$$
\end{corollary}

\begin{proof}
Note that multiplier ideals of monomial ideals are monomial
(for example, because we can find a log resolution which is
equivariant with respect to the standard $(\CC^*)^d$-action on $\AAA^d$).
If $Z$ is the subscheme defined by $T^{\bf u}$,
note that $T^{\bf u}\in\II(X,\alpha\cdot\fra)$
if and only if $(X,\lambda\cdot Y- Z)$ is log canonical
for some $\lambda>\alpha$. The assertion of the corollary
follows now from the above proposition.
\end{proof}

\providecommand{\bysame}{\leavevmode \hbox \o3em
{\hrulefill}\thinspace}

\end{document}